\documentclass[11pt]{amsart}
\usepackage{amsmath,amssymb,stmaryrd}
\newtheorem{theorem}{Theorem}
\newtheorem{proposition}[theorem]{Proposition}
\newtheorem{corollary}[theorem]{Corollary}
\newtheorem{lemma}[theorem]{Lemma}
\begin{document}

\title{Manifolds with $1/4$-pinched Curvature are Space Forms}
\author{Simon Brendle and Richard Schoen}
\address{Department of Mathematics \\
                 Stanford University \\
                 Stanford, CA 94305}
\thanks{The first author was partially supported by a Sloan Foundation Fellowship and by NSF grant DMS-0605223. The second author was partially supported by NSF grant DMS-0604960.}
\maketitle 


\section{Introduction}

One of the basic problems of Riemannian geometry is the classification of manifolds of
positive sectional curvature. The known examples include the spherical space forms
which carry constant curvature metrics and the rank $1$ symmetric spaces whose canonical
metrics have sectional curvatures at each point varying between $1$ and $4$. In 1951, 
H.E.~ Rauch \cite{Rauch} introduced the notion of curvature pinching for Riemannian manifolds
and posed the question of whether a compact, simply connected manifold $M$ whose sectional 
curvatures all lie in the interval $(1,4]$ is necessarily homeomorphic to the sphere $S^n$. 
This was proven by M.~Berger \cite{Berger} and W.~Klingenberg \cite{Klingenberg} around 1960 
using comparison techniques. However, this theorem leaves open the question of whether $M$ is diffeomorphic to $S^n$. This conjecture is known as the Differentiable Sphere Theorem, and the purpose of this paper is to prove this and a more general result which we describe.

We will say that a manifold $M$ has pointwise $1/4$-pinched sectional curvatures if $M$ has positive sectional curvature and for every point $p \in M$ the ratio of the maximum to the minimum sectional curvature at that point is less than $4$. In other words, for every pair of two-planes $\pi_1,\pi_2 \subset T_p M$ we have $0 < K(\pi_1) < 4 \, K(\pi_2)$. Our main result is the following:

\begin{theorem} 
\label{sphere.theorem}
Let $M$ be a compact Riemannian manifold of dimension $n \geq 4$ with pointwise $1/4$-pinched sectional curvatures. Then $M$ admits a metric of constant curvature and therefore is diffeomorphic to a spherical space form.
\end{theorem}

The techniques in this paper can be extended to give a classification of manifolds with weakly $1/4$-pinched sectional curvatures. We refer to \cite{Brendle-Schoen} for details.

Since our method of proof gives a canonical deformation from the $1/4$-pinched metric to a constant curvature metric, we can also prove the following equivariant version.

\begin{theorem}
\label{equivariant.sphere.theorem}
Let $M$ be a compact, simply connected Riemannian manifold of dimension $n \geq 4$ with pointwise $1/4$-pinched sectional curvatures. Assume that $G$ is a compact Lie group and $\rho$ is a group homomorphism from $G$ into the isometry group of $M$. Then there exists a group homomorphism $\sigma$ from $G$ into $O(n+1)$ and a diffeomorphism
$F$ from $M$ to $S^n$ which is equivariant; i.e. $F \circ \rho(g) = \sigma(g) \circ F$ for all $g \in G$.
\end{theorem}

Notice that Theorem \ref{sphere.theorem} and Theorem \ref{equivariant.sphere.theorem} do not assume any global pinching condition. A manifold is said to be globally $\delta$-pinched if all sectional curvatures at all points of $M$ lie in the interval $(1,\frac{1}{\delta}]$. The Differentiable Sphere Theorem under global $\delta$-pinching assumptions was obtained in 1966 by D.~Gromoll \cite{Gromoll} and E.~Calabi with a constant $\delta = \delta(n)$ converging to $1$ as $n \to \infty$. In 1971, M.~Sugimoto, K.~Shiohama, and H.~Karcher \cite{Sugimoto-Shiohama} proved the Differentiable Sphere Theorem with a pinching constant independent of $n$ ($\delta = 0.87$). The pinching constant was subsequently improved by E.~Ruh \cite{Ruh1} ($\delta = 0.80$) and by K.~Grove, H.~Karcher, and E.~Ruh \cite{Grove-Karcher-Ruh2} ($\delta = 0.76$). Ruh \cite{Ruh2} proved the Differentiable Sphere Theorem under a pointwise pinching condition, with a pinching constant converging to $1$ as $n \to \infty$.

The equivariant sphere theorem was first proven for globally $\delta$-pinched manifolds by K.~Grove, H.~Karcher, and E.~Ruh \cite{Grove-Karcher-Ruh1}, \cite{Grove-Karcher-Ruh2} with a pinching constant $\delta$ independent of $n$ ($\delta = 0.98$). The pinching constant was later improved by H.~Im~Hof and E.~Ruh \cite{Im-Hof-Ruh}. 

In 1982, R.~Hamilton \cite{Hamilton1} introduced a fundamental new tool to this problem. Given a compact Riemannian manifold $(M,g_0)$, Hamilton evolved the Riemannian metric by the equation 
\[\frac{\partial}{\partial t} g(t) = -2 \, \text{\rm Ric}_{g(t)}\] 
with initial condition $g(0) = g_0$. This equation is known as the Ricci flow. Hamilton also defined a normalized version of the Ricci flow. The normalized Ricci flow is defined by 
\[\frac{\partial}{\partial t} g(t) = -2 \, \text{\rm Ric}_{g(t)} + \frac{2}{n} \, r_{g(t)} \, g(t),\] 
where $r_{g(t)}$ denotes the mean value of the scalar curvature of $g(t)$. Note that the volume of $M$ is constant under the normalized flow. 

Using this method, Hamilton \cite{Hamilton1} proved that every three-manifold with positive Ricci curvature admits a constant curvature metric. In a subsequent paper, Hamilton \cite{Hamilton2} laid the general framework for the application of Ricci flow to Riemannian geometry and showed that four-manifolds with positive curvature operator are space forms. In 1991, H.~Chen \cite{Chen} extended this result to four-manifolds with $2$-positive curvature operator, which implies Theorem \ref{sphere.theorem} for $n = 4$. More recently, B.~Andrews and H.~Nguyen \cite{Andrews-Nguyen} proved that four-manifolds with $1/4$-pinched flag curvature are space forms. In higher dimensions, the Ricci flow was used by G.~Huisken \cite{Huisken} to show that sufficiently pinched manifolds are space forms (see also \cite{Margerin}, \cite{Nishikawa}). 

C.~B\"ohm and B.~Wilking \cite{Bohm-Wilking} used the Ricci flow to prove that manifolds with $2$-positive curvature operator are space forms. Most importantly, their work introduces new methods for deforming invariant sets and constructing pinching sets for the ODE on the space of curvature-type tensors arising from the evolution of the curvature. The curvature ODE was introduced by R. Hamilton \cite{Hamilton2} and had been exploited effectively in dimensions $3$ and $4$.

In 1988, M.~Micallef and J.D.~Moore \cite{Micallef-Moore} introduced minimal surface techniques into this problem and proved the topological sphere theorem for pointwise $1/4$-pinched manifolds using variational theory for the 
energy functional on maps from $S^2$ to $M$. Another important contribution of their paper was that they introduced a new curvature condition, positive isotropic curvature. This condition arose from consideration of the second variation of energy for maps of surfaces into $M$. The condition says that for every orthonormal four-frame $\{e_1,e_2,e_3,e_4\}$ we have the inequality
\[R_{1313} + R_{1414} + R_{2323} + R_{2424} - 2 \, R_{1234} > 0.\]
If we allow the weak inequality, then we say that $M$ has nonnegative isotropic curvature. Micallef and Moore proved that a compact, simply connected manifold with positive isotropic curvature is homeomorphic to $S^n$. Moreover, they observed
that pointwise $1/4$-pinching implies positive isotropic curvature. 

In dimension $4$, it was shown by R.~Hamilton \cite{Hamilton3} that positive isotropic curvature is preserved by the Ricci flow. The Ricci flow on four-manifolds with positive isotropic curvature will, in general, develop singularities. Hamilton established pointwise estimates for the curvature tensor of the evolving metric and used them to give a precise description of the singularities in this situation \cite{Hamilton-survey}, \cite{Hamilton3}. In order to extend the flow beyond singularities, Hamilton introduced the notion of Ricci flow with surgeries (see also \cite{Chen-Zhu},
\cite{Perelman1}, \cite{Perelman2}).

In Section 2 we prove that positive isotropic curvature is preserved by the Ricci flow in all dimensions. By the maximum principle (cf. \cite{Hamilton2}, Theorem 4.3), it suffices to show that positive isotropic curvature is preserved by the Hamilton ODE. We were not able to show that all isotropic curvatures improve under the ODE. Instead, we prove that the {\it minimum} isotropic curvature increases under the ODE, which is sufficient for our purposes. This is a very intricate calculation which exploits special identities and inequalities for the curvature tensor $R$ arising from the first and second variations applied to a set of four orthonormal vectors which minimize the isotropic curvature. After this paper was written, we learned that H.~Nguyen \cite{Nguyen} has independently proved that positive isotropic curvature is preserved under the Ricci flow. 

Even knowing that positive isotropic curvature is preserved, it seems to be a difficult analytic problem to give a complete
analysis of solutions to the Ricci flow satisfying that condition. A combination of results of M.~Micallef and M.~Wang \cite{Micallef-Wang} and recent results of A.~Fraser \cite{Fraser} on non-simply connected manifolds with positive isotropic curvature suggest that Hamilton's four dimensional flow with surgeries may hold in all dimensions. We do not treat this question here. Instead, we establish a convergence result for the Ricci flow in dimension $n \geq 4$ under a curvature condition which is substantially stronger than positive isotropic curvature, but includes the pointwise $1/4$-pinched manifolds. 

Given a curvature tensor $R$ thought of as a four-tensor on $\mathbb{R}^n$, we define $\tilde{R}$ to be the extension of $R$ as a curvature tensor on $\mathbb{R}^n \times \mathbb{R}$ which is zero in the additional direction. Thus $\tilde{R}$ is the curvature tensor one obtains for the manifold $M \times \mathbb{R}$. The condition that $\tilde{R}$ has positive isotropic curvature is preserved by the Ricci flow and is a much stronger condition than positive isotropic curvature itself. If $\tilde{R}$ has positive isotropic curvature, then $R$ has $2$-positive flag curvature in the sense that 
\[R_{1212} + R_{1313} > 0\] 
for all orthonormal three-frames $\{e_1,e_2,e_3\}$. In particular, this condition implies that $R$ has positive Ricci curvature. Continuing in this vein, we define $\hat{R}$ to be the curvature tensor on $\mathbb{R}^n \times \mathbb{R}^2$ obtained by extending $R$ to be zero in the two additional directions. Thus $\hat{R}$ is the curvature tensor of $M \times \mathbb{R}^2$. The condition that $\hat{R}$ has nonnegative isotropic curvature is, again, preserved by the Ricci flow. Moreover, this 
condition implies that $R$ has nonnegative sectional curvature. Note that $\hat{R}$ cannot have positive isotropic curvature due to the two flat directions.

This construction provides us with a convex cone in the space of algebraic curvature operators which is invariant under the Hamilton ODE, is contained in the cone of curvature operators with nonnegative sectional curvature, and contains all  nonnegative curvature operators. We may then directly apply results of B\"ohm and Wilking \cite{Bohm-Wilking} to obtain suitable pinching sets for the ODE.  Convergence of the normalized Ricci flow to a constant curvature metric then follows from work of
Hamilton \cite{Hamilton2} (see also \cite{Bohm-Wilking}). This material is discussed in detail in Section 3.

Finally, in Section 4 we give a necessary and sufficient condition for $\hat{R}$ to have nonnegative isotropic curvature. This yields the following result:

\begin{theorem} 
\label{convergence.of.ricci.flow}
Let $(M,g_0)$ be a compact Riemannian manifold of dimension $n \geq 4$. Assume that 
\[R_{1313} + \lambda^2 \, R_{1414} + \mu^2 \, R_{2323} + \lambda^2\mu^2 \, R_{2424} - 2\lambda\mu \, R_{1234} > 0\] 
for all orthonormal four-frames $\{e_1,e_2,e_3,e_4\}$ and all $\lambda,\mu \in [-1,1]$. Then the normalized 
Ricci flow with initial metric $g_0$ exists for all time and converges to a constant curvature metric as $t \to \infty$.
\end{theorem}

It follows from Berger's inequality that every manifold with pointwise $1/4$-pinched sectional curvatures 
satisfies the curvature condition in Theorem \ref{convergence.of.ricci.flow}. Hence, Theorem \ref{sphere.theorem} and Theorem \ref{equivariant.sphere.theorem} are immediate consequences of Theorem \ref{convergence.of.ricci.flow}.

\section{Positive isotropic curvature is preserved by the Ricci flow}

In this section, we will prove that positive isotropic curvature is preserved by the Ricci flow. By work of R.~Hamilton \cite{Hamilton2}, it suffices to show that positive isotropic curvature is preserved by the ODE $\frac{d}{dt} R = Q(R)$, where $Q(R)$ is defined by 
\[Q(R)_{ijkl} = R_{ijpq} \, R_{klpq} + 2 \, R_{ipkq} \, R_{jplq} - 2 \, R_{iplq} \, R_{jpkq}.\] 
To that end, we assume that $R$ is a curvature tensor with nonnegative isotropic curvature. Moreover, suppose that $\{e_1,e_2,e_3,e_4\}$ is an orthonormal four-frame satisfying 
\begin{equation} 
\label{zero.isotropic.curvature}
R_{1313} + R_{1414} + R_{2323} + R_{2424} - 2 \, R_{1234} = 0. 
\end{equation}
We will show that 
\begin{equation} 
Q(R)_{1313} + Q(R)_{1414} + Q(R)_{2323} + Q(R)_{2424} - 2 \, Q(R)_{1234} \geq 0. 
\end{equation}
The following observation will be useful: if $\{e_1,e_2,e_3,e_4\}$ is an orthonormal four-frame satisfying (\ref{zero.isotropic.curvature}), then the four-frames $\{e_2,-e_1,e_3,e_4\}$, $\{e_2,-e_1,e_4,-e_3\}$, and $\{e_3,e_4,e_1,e_2\}$ also satisfy (\ref{zero.isotropic.curvature}). Hence, any statement that we can prove for the frame $\{e_1,e_2,e_3,e_4\}$ will also hold for the frames $\{e_2,-e_1,e_3,e_4\}$, $\{e_2,-e_1,e_4,-e_3\}$, and $\{e_3,e_4,e_1,e_2\}$.

Using the first Bianchi identity, we obtain 
\begin{align*} 
Q(R)_{1234} 
&= R_{12pq} \, R_{34pq} + 2 \, R_{1p3q} \, R_{2p4q} - 2 \, R_{1p4q} \, R_{2p3q} \\ 
&= R_{12pq} \, R_{34pq} + R_{13pq} \, R_{24pq} - R_{14pq} \, R_{23pq} \\ 
&+ 2 \, R_{1p3q} \, R_{4p2q} - 2 \, R_{1p4q} \, R_{3p2q}. 
\end{align*}
This implies 
\begin{align*} 
&Q(R)_{1313} + Q(R)_{1414} + Q(R)_{2323} + Q(R)_{2424} - 2 \, Q(R)_{1234} \\ 
&= R_{13pq} \, R_{13pq} + 2 \, R_{1p1q} \, R_{3p3q} - 2 \, R_{1p3q} \, R_{3p1q} \\ 
&+ R_{14pq} \, R_{14pq} + 2 \, R_{1p1q} \, R_{4p4q} - 2 \, R_{1p4q} \, R_{4p1q} \\ 
&+ R_{23pq} \, R_{23pq} + 2 \, R_{2p2q} \, R_{3p3q} - 2 \, R_{2p3q} \, R_{3p2q} \\ 
&+ R_{24pq} \, R_{24pq} + 2 \, R_{2p2q} \, R_{4p4q} - 2 \, R_{2p4q} \, R_{4p2q} \\ 
&- 2 \, R_{12pq} \, R_{34pq} - 2 \, R_{13pq} \, R_{24pq} + 2 \, R_{14pq} \, R_{23pq} \\ 
&- 4 \, R_{1p3q} \, R_{4p2q} + 4 \, R_{1p4q} \, R_{3p2q}. 
\end{align*} 
Rearranging terms yields 
\begin{align*} 
&Q(R)_{1313} + Q(R)_{1414} + Q(R)_{2323} + Q(R)_{2424} - 2 \, Q(R)_{1234} \\ 
&= (R_{13pq} - R_{24pq}) \, (R_{13pq} - R_{24pq}) + (R_{14pq} + R_{23pq}) \, (R_{14pq} + R_{23pq}) \\ 
&+ 2 \, (R_{1p1q} + R_{2p2q}) \, (R_{3p3q} + R_{4p4q}) - 2 \, R_{12pq} \, R_{34pq} \\ 
&- 2 \, (R_{1p3q} + R_{2p4q}) \, (R_{3p1q} + R_{4p2q}) - 2 \, (R_{1p4q} - R_{2p3q}) \, (R_{4p1q} - R_{3p2q}). 
\end{align*} 
The first two terms on the right are clearly nonnegative. 

\begin{lemma}
\label{first.order.conditions.1}
We have 
\[R_{1213} + R_{1242} + R_{3413} + R_{3442} = R_{1214} + R_{1223} + R_{3414} + R_{3423} = 0.\]
\end{lemma} 

\textbf{Proof.} 
Consider the frame $\{e_1,\cos s \, e_2 - \sin s \, e_3,\sin s \, e_2 + \cos s \, e_3,e_4\}$. Since $R$ has nonnegative isotropic curvature, the function 
\begin{align*} 
s \mapsto &\cos^2s \, (R_{1313} + R_{2424} - 2 \, R_{1234}) + \sin^2s \, (R_{1212} + R_{3434} + 2 \, R_{1324}) \\ &+ R_{1414} + R_{2323} + 2\cos s \, \sin s \, (R_{1213} - R_{2434} - R_{1224} + R_{1334}) 
\end{align*} 
is nonnegative and vanishes for $s = 0$. This implies $R_{1213} - R_{2434} - R_{1224} + R_{1334} = 0$. If we replace $\{e_1,e_2,e_3,e_4\}$ by $\{e_2,-e_1,e_3,e_4\}$, we obtain 
$-R_{2123} + R_{1434} - R_{2114} + R_{2334} = 0$. \\

\begin{proposition} 
\label{p.and.q.between.1.and.4}
We have 
\begin{align*} 
&\sum_{p,q=1}^4 (R_{1p1q} + R_{2p2q}) \, (R_{3p3q} + R_{4p4q}) - \sum_{p,q=1}^4 R_{12pq} \, R_{34pq} \\ 
&= \sum_{p,q=1}^4 (R_{1p3q} + R_{2p4q}) \, (R_{3p1q} + R_{4p2q}) \\ &+ \sum_{p,q=1}^4 (R_{1p4q} - R_{2p3q}) \, (R_{4p1q} - R_{3p2q}). 
\end{align*} 
\end{proposition}

\textbf{Proof.} 
Direct computation yields 
\begin{align*} 
&\sum_{p,q=1}^4 (R_{1p1q} + R_{2p2q}) \, (R_{3p3q} + R_{4p4q}) - \sum_{p,q=1}^4 R_{12pq} \, R_{34pq} \\ 
&- \sum_{p,q=1}^4 (R_{1p3q} + R_{2p4q}) \, (R_{3p1q} + R_{4p2q}) \\ &- \sum_{p,q=1}^4 (R_{1p4q} - R_{2p3q}) \, (R_{4p1q} - R_{3p2q}) \\ 
&= (R_{1212} + R_{3434}) \, (R_{1313} + R_{1414} + R_{2323} + R_{2424} - 2 \, R_{1234}) \\ 
&+ 2 \, R_{1234} \, (R_{1313} + R_{1414} + R_{2323} + R_{2424} + 2 \, R_{1342} + 2 \, R_{1423}) \\ 
&- (R_{1213} + R_{1242} + R_{3413} + R_{3442})^2 - (R_{1214} + R_{1223} + R_{3414} + R_{3423})^2 \\
&= (R_{1212} + R_{3434} + 2 \, R_{1234}) \, (R_{1313} + R_{1414} + R_{2323} + R_{2424} - 2 \, R_{1234}) \\ 
&- (R_{1213} + R_{1242} + R_{3413} + R_{3442})^2 - (R_{1214} + R_{1223} + R_{3414} + R_{3423})^2. 
\end{align*} 
The expression on the right is zero by Lemma \ref{first.order.conditions.1}. \\

\begin{lemma} 
\label{first.order.conditions.2}
We have 
\[R_{133q} + R_{144q} + R_{432q} = R_{233q} + R_{244q} + R_{341q} = 0\] 
for all $5 \leq q \leq n$.
\end{lemma}

\textbf{Proof.} 
Consider the frame $\{\cos s \, e_1 + \sin s \, e_q,e_2,e_3,e_4\}$. Since $R$ has nonnegative isotropic curvature, the function 
\begin{align*} 
s \mapsto &\cos^2s \, (R_{1313} + R_{1414}) + \sin^2s \, (R_{q3q3} + R_{q4q4}) + R_{2323} + R_{2424} \\ 
&+ 2\cos s \, \sin s \, (R_{13q3} + R_{14q4}) - 2\cos s \, R_{1234} - 2\sin s \, R_{q234} 
\end{align*} 
is nonnegative, and vanishes for $s = 0$. This implies $R_{13q3} + R_{14q4} - R_{q234} = 0$. If we replace $\{e_1,e_2,e_3,e_4\}$ by $\{e_2,-e_1,e_3,e_4\}$, we obtain $R_{23q3} + R_{24q4} + R_{q134} = 0$. \\

\begin{proposition} 
\label{p.between.1.and.4.q.bigger.than.4}
Fix $q$ such that $5 \leq q \leq n$. Then we have 
\begin{align*} 
&\sum_{p=1}^4 (R_{1p1q} + R_{2p2q}) \, (R_{3p3q} + R_{4p4q}) - \sum_{p=1}^4 R_{12pq} \, R_{34pq} \\ 
&= \sum_{p=1}^4 (R_{1p3q} + R_{2p4q}) \, (R_{3p1q} + R_{4p2q}) \\ &+ \sum_{p=1}^4 (R_{1p4q} - R_{2p3q}) \, (R_{4p1q} - R_{3p2q}). 
\end{align*} 
\end{proposition}

\textbf{Proof.} 
Using Lemma \ref{first.order.conditions.2}, we obtain 
\begin{align*} 
&\sum_{p=1}^2 (R_{1p1q} + R_{2p2q}) \, (R_{3p3q} + R_{4p4q}) - \sum_{p=1}^2 R_{12pq} \, R_{34pq} \\ 
&= R_{212q} \, (R_{313q} + R_{414q}) + R_{121q} \, (R_{323q} + R_{424q}) \\ &- R_{121q} \, R_{341q} - R_{122q} \, R_{342q} \\ 
&= R_{212q} \, (R_{313q} + R_{414q} + R_{342q}) \\ &+ R_{121q} \, (R_{323q} + R_{424q} - R_{341q}) \\ 
&= 0
\end{align*} 
and 
\begin{align*} 
&\sum_{p=3}^4 (R_{1p3q} + R_{2p4q}) \, (R_{3p1q} + R_{4p2q}) \\ &+ \sum_{p=3}^4 (R_{1p4q} - R_{2p3q}) \, (R_{4p1q} - R_{3p2q}) \\ 
&= (R_{133q} + R_{234q}) \, R_{432q} + (R_{143q} + R_{244q}) \, R_{341q} \\ 
&+ (R_{134q} - R_{233q}) \, R_{431q} - (R_{144q} - R_{243q}) \, R_{342q} \\ 
&= (R_{133q} + R_{234q} + R_{144q} - R_{243q}) \, R_{432q} \\ 
&+ (R_{143q} + R_{244q} - R_{134q} + R_{233q}) \, R_{341q} \\ 
&= (R_{133q} + R_{144q} + R_{432q}) \, R_{432q} \\ 
&+ (R_{341q} + R_{244q} + R_{233q}) \, R_{341q} \\ 
&= 0. 
\end{align*}
Replacing $\{e_1,e_2,e_3,e_4\}$ by $\{e_3,e_4,e_1,e_2\}$ yields 
\[\sum_{p=3}^4 (R_{1p1q} + R_{2p2q}) \, (R_{3p3q} + R_{4p4q}) - \sum_{p=3}^4 R_{12pq} \, R_{34pq} = 0\] 
and 
\begin{align*} 
&\sum_{p=1}^2 (R_{1p3q} + R_{2p4q}) \, (R_{3p1q} + R_{4p2q}) \\ 
&+ \sum_{p=1}^2 (R_{1p4q} - R_{2p3q}) \, (R_{4p1q} - R_{3p2q}) = 0. 
\end{align*} 
Putting these facts together, the assertion follows. \\

\begin{proposition} 
\label{min}
Assume that $w_1,w_2,w_3,w_4$ are orthogonal to $e_1,e_2,e_3,e_4$. Then the expression 
\begin{align*} 
&R(w_1,e_3,w_1,e_3) + R(w_1,e_4,w_1,e_4) \\ 
&+ R(w_2,e_3,w_2,e_3) + R(w_2,e_4,w_2,e_4) \\ 
&+ R(e_1,w_3,e_1,w_3) + R(e_2,w_3,e_2,w_3) \\ 
&+ R(e_1,w_4,e_1,w_4) + R(e_2,w_4,e_2,w_4) \\ 
&- 2 \, \big [ R(e_3,w_1,e_1,w_3) + R(e_4,w_1,e_2,w_3) \big ] \\ 
&- 2 \, \big [ R(e_4,w_1,e_1,w_4) - R(e_3,w_1,e_2,w_4) \big ] \\ 
&+ 2 \, \big [ R(e_4,w_2,e_1,w_3) - R(e_3,w_2,e_2,w_3) \big ] \\ 
&- 2 \, \big [ R(e_3,w_2,e_1,w_4) + R(e_4,w_2,e_2,w_4) \big ] \\ 
&- 2 \, R(w_1,w_2,e_3,e_4) - 2 \, R(e_1,e_2,w_3,w_4) 
\end{align*} 
is nonnegative.
\end{proposition}

\textbf{Proof.} 
For $i = 1, \hdots, 4$, we denote by $v_i(s)$ the solution of the ODE 
\[v_i'(s) = \sum_{j=1}^4 (\langle v_i(s),e_j \rangle \, w_j - \langle v_i(s),w_j \rangle \, e_j)\] 
with initial condition $v_i(0) = e_i$. Clearly, $\{v_1(s),v_2(s),v_3(s),v_4(s)\}$ is an orthonormal four-frame. Moreover, 
$v_i'(0) = w_i$ and $v_i''(0) = -\sum_{j=1}^4 \langle w_i,w_j \rangle \, e_j$. Since $R$ has nonnegative isotropic curvature, the function 
\begin{align*} 
s \mapsto &R(v_1(s),v_3(s),v_1(s),v_3(s)) + R(v_1(s),v_4(s),v_1(s),v_4(s)) \\ 
&+ R(v_2(s),v_3(s),v_2(s),v_3(s)) + R(v_2(s),v_4(s),v_2(s),v_4(s)) \\ 
&- 2 \, R(v_1(s),v_2(s),v_3(s),v_4(s)) 
\end{align*} 
is nonnegative and vanishes for $s = 0$. Therefore, the second derivative of this function at $s = 0$ is nonnegative. This implies 
\[0 \leq J^{(1)} + J^{(2)} + J^{(3)} + J^{(4)} - J^{(5)},\] 
where 
\begin{align*} 
J^{(1)} &= \frac{1}{2} \, \frac{d^2}{ds^2} R(v_1(s),v_3(s),v_1(s),v_3(s)) \Big |_{s=0} \\ 
&= R(w_1,e_3,w_1,e_3) + R(e_1,w_3,e_1,w_3) \\ 
&+ 2 \, R(e_1,e_3,w_1,w_3) + 2 \, R(e_1,w_3,w_1,e_3) \\ 
&- (|w_1|^2 + |w_3|^2) \, R(e_1,e_3,e_1,e_3) \\ 
&- \langle w_1,w_2 \rangle \, R(e_1,e_3,e_2,e_3) - \langle w_1,w_4 \rangle \, R(e_1,e_3,e_4,e_3) \\ 
&- \langle w_3,w_2 \rangle \, R(e_1,e_3,e_1,e_2) - \langle w_3,w_4 \rangle \, R(e_1,e_3,e_1,e_4), 
\end{align*} 
\begin{align*} 
J^{(2)} &= \frac{1}{2} \, \frac{d^2}{ds^2} R(v_1(s),v_4(s),v_1(s),v_4(s)) \Big |_{s=0} \\ 
&= R(w_1,e_4,w_1,e_4) + R(e_1,w_4,e_1,w_4) \\ 
&+ 2 \, R(e_1,e_4,w_1,w_4) + 2 \, R(e_1,w_4,w_1,e_4) \\ 
&- (|w_1|^2 + |w_4|^2) \, R(e_1,e_4,e_1,e_4) \\ 
&- \langle w_1,w_2 \rangle \, R(e_1,e_4,e_2,e_4) - \langle w_1,w_3 \rangle \, R(e_1,e_4,e_3,e_4) \\ 
&- \langle w_4,w_2 \rangle \, R(e_1,e_4,e_1,e_2) - \langle w_4,w_3 \rangle \, R(e_1,e_4,e_1,e_3), 
\end{align*} 
\begin{align*} 
J^{(3)} &= \frac{1}{2} \, \frac{d^2}{ds^2} R(v_2(s),v_3(s),v_2(s),v_3(s)) \Big |_{s=0} \\ 
&= R(w_2,e_3,w_2,e_3) + R(e_2,w_3,e_2,w_3) \\ 
&+ 2 \, R(e_2,e_3,w_2,w_3) + 2 \, R(e_2,w_3,w_2,e_3) \\ 
&- (|w_2|^2 + |w_3|^2) \, R(e_2,e_3,e_2,e_3) \\ 
&- \langle w_2,w_1 \rangle \, R(e_2,e_3,e_1,e_3) - \langle w_2,w_4 \rangle \, R(e_2,e_3,e_4,e_3) \\ 
&- \langle w_3,w_1 \rangle \, R(e_2,e_3,e_2,e_1) - \langle w_3,w_4 \rangle \, R(e_2,e_3,e_2,e_4), 
\end{align*} 
\begin{align*} 
J^{(4)} &= \frac{1}{2} \, \frac{d^2}{ds^2} R(v_2(s),v_4(s),v_2(s),v_4(s)) \Big |_{s=0} \\ 
&= R(w_2,e_4,w_2,e_4) + R(e_2,w_4,e_2,w_4) \\ 
&+ 2 \, R(e_2,e_4,w_2,w_4) + 2 \, R(e_2,w_4,w_2,e_4) \\ 
&- (|w_2|^2 + |w_4|^2) \, R(e_2,e_4,e_2,e_4) \\ 
&- \langle w_2,w_1 \rangle \, R(e_2,e_4,e_1,e_4) - \langle w_2,w_3 \rangle \, R(e_2,e_4,e_3,e_4) \\ 
&- \langle w_4,w_1 \rangle \, R(e_2,e_4,e_2,e_1) - \langle w_4,w_3 \rangle \, R(e_2,e_4,e_2,e_3), 
\end{align*} 
and 
\begin{align*} 
J^{(5)} &= \frac{d^2}{ds^2} R(v_1(s),v_2(s),v_3(s),v_4(s)) \Big |_{s=0} \\ 
&= 2 \, R(w_1,w_2,e_3,e_4) + 2 \, R(w_1,e_2,w_3,e_4) + 2 \, R(w_1,e_2,e_3,w_4) \\ 
&+ 2 \, R(e_1,w_2,w_3,e_4) + 2 \, R(e_1,w_2,e_3,w_4) + 2 \, R(e_1,e_2,w_3,w_4) \\ 
&- (|w_1|^2 + |w_2|^2 + |w_3|^2 + |w_4|^2) \, R(e_1,e_2,e_3,e_4) \\ 
&- \langle w_1,w_3 \rangle \, R(e_3,e_2,e_3,e_4) - \langle w_1,w_4 \rangle \, R(e_4,e_2,e_3,e_4) \\ 
&- \langle w_2,w_3 \rangle \, R(e_1,e_3,e_3,e_4) - \langle w_2,w_4 \rangle \, R(e_1,e_4,e_3,e_4) \\ 
&- \langle w_3,w_1 \rangle \, R(e_1,e_2,e_1,e_4) - \langle w_3,w_2 \rangle \, R(e_1,e_2,e_2,e_4) \\ 
&- \langle w_4,w_1 \rangle \, R(e_1,e_2,e_3,e_1) - \langle w_4,w_2 \rangle \, R(e_1,e_2,e_3,e_2). 
\end{align*} 
Rearranging terms yields 
\begin{align*} 
0 &\leq R(w_1,e_3,w_1,e_3) + R(w_1,e_4,w_1,e_4) \\ 
&+ R(w_2,e_3,w_2,e_3) + R(w_2,e_4,w_2,e_4) \\ 
&+ R(e_1,w_3,e_1,w_3) + R(e_2,w_3,e_2,w_3) \\ 
&+ R(e_1,w_4,e_1,w_4) + R(e_2,w_4,e_2,w_4) \\ 
&+ 2 \, R(e_1,e_3,w_1,w_3) + 2 \, R(e_1,w_3,w_1,e_3) - 2 \, R(w_1,e_2,w_3,e_4) \\ 
&+ 2 \, R(e_1,e_4,w_1,w_4) + 2 \, R(e_1,w_4,w_1,e_4) - 2 \, R(w_1,e_2,e_3,w_4) \\ 
&+ 2 \, R(e_2,e_3,w_2,w_3) + 2 \, R(e_2,w_3,w_2,e_3) - 2 \, R(e_1,w_2,w_3,e_4) \\ 
&+ 2 \, R(e_2,e_4,w_2,w_4) + 2 \, R(e_2,w_4,w_2,e_4) - 2 \, R(e_1,w_2,e_3,w_4) \\ 
&- 2 \, R(w_1,w_2,e_3,e_4) - 2 \, R(e_1,e_2,w_3,w_4) \\ 
&- |w_1|^2 \, (R_{1313} + R_{1414} - R_{1234}) - |w_2|^2 \, (R_{2323} + R_{2424} - R_{1234}) \\ 
&- |w_3|^2 \, (R_{1313} + R_{2323} - R_{1234}) - |w_4|^2 \, (R_{1414} + R_{2424} - R_{1234}) \\ 
&+ (\langle w_1,w_3 \rangle - \langle w_2,w_4 \rangle) \, (R_{1214} - R_{1232} + R_{3234} - R_{1434}) \\ 
&- (\langle w_1,w_4 \rangle + \langle w_2,w_3 \rangle) \, (R_{1213} + R_{1242} + R_{3134} + R_{2434}) \\ 
&- 2 \, \langle w_1,w_2 \rangle \, (R_{1323} + R_{1424}) - 2 \, \langle w_3,w_4 \rangle \, (R_{1314} + R_{2324}). 
\end{align*}  
We now replace the frame $\{e_1,e_2,e_3,e_4\}$ by $\{e_2,-e_1,e_4,-e_3\}$. This yields 
\begin{align*} 
0 &\leq R(w_1,e_4,w_1,e_4) + R(w_1,e_3,w_1,e_3) \\ 
&+ R(w_2,e_4,w_2,e_4) + R(w_2,e_3,w_2,e_3) \\ 
&+ R(e_2,w_3,e_2,w_3) + R(e_1,w_3,e_1,w_3) \\ 
&+ R(e_2,w_4,e_2,w_4) + R(e_1,w_4,e_1,w_4) \\ 
&+ 2 \, R(e_2,e_4,w_1,w_3) + 2 \, R(e_2,w_3,w_1,e_4) - 2 \, R(w_1,e_1,w_3,e_3) \\ 
&- 2 \, R(e_2,e_3,w_1,w_4) - 2 \, R(e_2,w_4,w_1,e_3) + 2 \, R(w_1,e_1,e_4,w_4) \\ 
&- 2 \, R(e_1,e_4,w_2,w_3) - 2 \, R(e_1,w_3,w_2,e_4) + 2 \, R(e_2,w_2,w_3,e_3) \\ 
&+ 2 \, R(e_1,e_3,w_2,w_4) + 2 \, R(e_1,w_4,w_2,e_3) - 2 \, R(e_2,w_2,e_4,w_4) \\ 
&+ 2 \, R(w_1,w_2,e_4,e_3) + 2 \, R(e_2,e_1,w_3,w_4) \\ 
&- |w_1|^2 \, (R_{2424} + R_{2323} - R_{2143}) - |w_2|^2 \, (R_{1414} + R_{1313} - R_{2143}) \\ 
&- |w_3|^2 \, (R_{2424} + R_{1414} - R_{2143}) - |w_4|^2 \, (R_{2323} + R_{1313} - R_{2143}) \\ 
&+ (\langle w_1,w_3 \rangle - \langle w_2,w_4 \rangle) \, (R_{2123} - R_{2141} + R_{4143} - R_{2343}) \\ 
&+ (\langle w_1,w_4 \rangle + \langle w_2,w_3 \rangle) \, (R_{2124} + R_{2131} + R_{4243} + R_{1343}) \\ 
&+ 2 \, \langle w_1,w_2 \rangle \, (R_{2414} + R_{2313}) + 2 \, \langle w_3,w_4 \rangle \, (R_{2423} + R_{1413}). 
\end{align*}  
In the next step, we take the arithmetic mean of both inequalitities. Using the identity $R_{1313} + R_{1414} + R_{2323} + R_{2424} - 2 \, R_{1234} = 0$, we obtain 
\begin{align*} 
0 &\leq R(w_1,e_3,w_1,e_3) + R(w_1,e_4,w_1,e_4) \\ 
&+ R(w_2,e_3,w_2,e_3) + R(w_2,e_4,w_2,e_4) \\ 
&+ R(e_1,w_3,e_1,w_3) + R(e_2,w_3,e_2,w_3) \\ 
&+ R(e_1,w_4,e_1,w_4) + R(e_2,w_4,e_2,w_4) \\ 
&+ \big [ R(e_1,e_3,w_1,w_3) + R(e_1,w_3,w_1,e_3) - R(w_1,e_2,w_3,e_4) \\ 
&\hspace{2mm} + R(e_2,e_4,w_1,w_3) + R(e_2,w_3,w_1,e_4) - R(w_1,e_1,w_3,e_3) \big ] \\ 
&+ \big [ R(e_1,e_4,w_1,w_4) + R(e_1,w_4,w_1,e_4) - R(w_1,e_2,e_3,w_4) \\ 
&\hspace{2mm} - R(e_2,e_3,w_1,w_4) - R(e_2,w_4,w_1,e_3) + R(w_1,e_1,e_4,w_4) \big ] \\ 
&+ \big [ R(e_2,e_3,w_2,w_3) + R(e_2,w_3,w_2,e_3) - R(e_1,w_2,w_3,e_4) \\ 
&\hspace{2mm} - R(e_1,e_4,w_2,w_3) - R(e_1,w_3,w_2,e_4) + R(e_2,w_2,w_3,e_3) \big ] \\ 
&+ \big [ R(e_2,e_4,w_2,w_4) + R(e_2,w_4,w_2,e_4) - R(e_1,w_2,e_3,w_4) \\ 
&\hspace{2mm} + R(e_1,e_3,w_2,w_4) + R(e_1,w_4,w_2,e_3) - R(e_2,w_2,e_4,w_4) \big ] \\ 
&- 2 \, R(w_1,w_2,e_3,e_4) - 2 \, R(e_1,e_2,w_3,w_4). 
\end{align*} 
The assertion follows now from the first Bianchi identity. \\

\begin{proposition} 
\label{p.and.q.bigger.than.4}
We have 
\begin{align*} 
&\sum_{p,q=5}^n (R_{1p1q} + R_{2p2q}) \, (R_{3p3q} + R_{4p4q}) - \sum_{p,q=5}^n R_{12pq} \, R_{34pq} \\ 
&\geq \sum_{p,q=5}^n (R_{1p3q} + R_{2p4q}) \, (R_{3p1q} + R_{4p2q}) \\ 
&+ \sum_{p,q=5}^n (R_{1p4q} - R_{2p3q}) \, (R_{4p1q} - R_{3p2q}). 
\end{align*}
\end{proposition}

\textbf{Proof.} 
Consider the following $(n-4) \times (n-4)$ matrices: 
\[\begin{array}{l@{\qquad}l} 
a_{pq} = R_{1p1q} + R_{2p2q}, & b_{pq} = R_{3p3q} +  R_{4p4q}, \\ 
c_{pq} = R_{3p1q} + R_{4p2q}, & d_{pq} = R_{4p1q} - R_{3p2q}, \\ 
e_{pq} = R_{12pq}, & f_{pq} = R_{34pq} 
\end{array}\] 
($5 \leq p,q \leq n$). It follows from Proposition \ref{min} that the matrix 
\[L = \begin{bmatrix} B & -F & -C & -D \\ F & B & D & -C \\ -C^T & D^T & A & -E \\ -D^T & -C^T & E & A \end{bmatrix}\] 
is positive semi-definite. We next define 
\[U = \begin{bmatrix} 0 & 0 & I & 0 \\ 0 & 0 & 0 & -I \\ -I & 0 & 0 & 0 \\ 0 & I & 0 & 0 \end{bmatrix}.\] 
Since $L$ is positive semi-definite, we have 
\begin{align*} 
0 &\leq \frac{1}{4} \, \text{\rm tr}(LULU^T) \\ 
&= \text{\rm tr}(AB) + \text{\rm tr}(EF) - \text{\rm tr}(C^2) - \text{\rm tr}(D^2) \\ 
&= \sum_{p,q=5}^n a_{pq} \, b_{pq} - \sum_{p,q=5}^n e_{pq} \, f_{pq} - \sum_{p,q=5}^n c_{pq} \, c_{qp} - \sum_{p,q=5}^n d_{pq} \, d_{qp}. 
\end{align*} 
This completes the proof. \\

\begin{corollary} 
\label{pic.is.preserved.1}
We have 
\begin{align*} 
&\sum_{p,q=1}^n (R_{1p1q} + R_{2p2q}) \, (R_{3p3q} + R_{4p4q}) - \sum_{p,q=1}^n R_{12pq} \, R_{34pq} \\ 
&\geq \sum_{p,q=1}^n (R_{1p3q} + R_{2p4q}) \, (R_{3p1q} + R_{4p2q}) \\ &+ \sum_{p,q=1}^n (R_{1p4q} - R_{2p3q}) \, (R_{4p1q} - R_{3p2q}). 
\end{align*} 
Consequently, 
\[Q(R)_{1313} + Q(R)_{1414} + Q(R)_{2323} + Q(R)_{2424} - 2 \, Q(R)_{1234} \geq 0.\] 
\end{corollary}

After these preparations, we now prove that nonnegative isotropic curvature is preserved by the ODE $\frac{d}{dt} R = Q(R)$:

\begin{proposition} 
\label{pic.is.preserved.2} 
Suppose that $R(t)$, $t \in [0,T)$, is a solution of the ODE $\frac{d}{dt} R(t) = Q(R(t))$. If $R(0)$ has nonnegative isotropic curvature, then $R(t)$ has nonnegative isotropic curvature for all $t \in [0,T)$.
\end{proposition}

\textbf{Proof.} 
Fix $\varepsilon > 0$, and denote by $R_\varepsilon(t)$ the solution of the ODE $\frac{d}{dt} R_\varepsilon(t) = Q(R_\varepsilon(t)) + \varepsilon I$ with initial condition $R_\varepsilon(0) = R(0) + \varepsilon I$. The function $R_\varepsilon(t)$ is defined on some time interval $[0,T_\varepsilon)$. We claim that $R_\varepsilon(t)$ has positive isotropic for all $t \in [0,T_\varepsilon)$. To prove this, we argue by contradiction. Suppose that there exists a time $t \in [0,T_\varepsilon)$ such that $R_\varepsilon(t)$ does not have positive isotropic curvature. Let 
\[\tau = \inf \{t \in [0,T_\varepsilon): \text{$R_\varepsilon(t)$ does not have positive isotropic curvature}\}.\] 
Clearly, $\tau > 0$. Moreover, there exists an orthonormal four-frame $\{e_1,e_2,e_3,e_4\}$ such that 
\[R_\varepsilon(\tau)_{1313} + R_\varepsilon(\tau)_{1414} + R_\varepsilon(\tau)_{2323} + R_\varepsilon(\tau)_{2424} - 2 \, R_\varepsilon(\tau)_{1234} = 0.\] 
By definition of $\tau$, $R_\varepsilon(t)$ has positive isotropic curvature for all $t \in [0,\tau)$. This implies 
\[R_\varepsilon(\tau)_{1313} + R_\varepsilon(\tau)_{1414} + R_\varepsilon(\tau)_{2323} + R_\varepsilon(\tau)_{2424} - 2 \, R_\varepsilon(\tau)_{1234} > 0\] 
for all $t \in [0,\tau)$. Therefore, we obtain 
\begin{align*} 
&Q(R_\varepsilon(\tau))_{1313} + Q(R_\varepsilon(\tau))_{1414} \\ 
&+ Q(R_\varepsilon(\tau))_{2323} + Q(R_\varepsilon(\tau))_{2424} - 2 \, Q(R_\varepsilon(\tau))_{1234} + 4\varepsilon \leq 0. 
\end{align*} 
On the other hand, since $R_\varepsilon(\tau)$ has nonnegative isotropic curvature, we have 
\begin{align*} 
&Q(R_\varepsilon(\tau))_{1313} + Q(R_\varepsilon(\tau))_{1414} \\ 
&+ Q(R_\varepsilon(\tau))_{2323} + Q(R_\varepsilon(\tau))_{2424} - 2 \, Q(R_\varepsilon(\tau))_{1234} \geq 0 
\end{align*} 
by Corollary \ref{pic.is.preserved.1}. This is a contradiction. \\

Therefore, $R_\varepsilon(t)$ has positive isotropic curvature for all $t \in [0,T_\varepsilon)$. Standard ODE theory implies that $T \leq \liminf_{\varepsilon \to 0} T_\varepsilon$ and $R(t) = \lim_{\varepsilon \to 0} R_\varepsilon(t)$ for all $t \in [0,T)$. Consequently, $R(t)$ has nonnegative isotropic curvature for all $t \in [0,T)$. \\

\section{Another invariant curvature condition for the Ricci flow} 

In this section, we construct a continuous family of cones that serves as a pinching family. Given any algebraic curvature operator $R$ on $\mathbb{R}^n$, we define an algebraic curvature operator $\hat{R}$ on $\mathbb{R}^n \times \mathbb{R}^2$ by 
\[\hat{R}(\hat{v}_1,\hat{v}_2,\hat{v}_3,\hat{v}_4) = R(v_1,v_2,v_3,v_4)\] 
for all vectors $\hat{v}_j = (v_j,x_j) \in \mathbb{R}^n \times \mathbb{R}^2$. We denote by $\hat{C}$ the set of all algebraic curvature operators on $\mathbb{R}^n$ with the property that $\hat{R}$ has nonnegative isotropic curvature: 
\[\hat{C} = \{R \in S_B^2(\mathfrak{so}(n)): \text{$\hat{R}$ has nonnegative isotropic curvature}\}.\] 
Clearly, $\hat{C}$ is a closed, convex, $O(n)$-invariant cone in the space of algebraic curvature operators. We next establish some basic properties of the cone $\hat{C}$:

\begin{proposition} 
\label{properties.of.C.hat}
The cone $\hat{C}$ has the following properties:
\begin{itemize}
\item[(i)] The cone $\hat{C}$ is invariant under the ODE $\frac{d}{dt} R = Q(R)$.
\item[(ii)] Every algebraic curvature operator $R \in \hat{C}$ has nonnegative sectional curvature.
\item[(iii)] If $R$ is a nonnegative curvature operator on $\mathbb{R}^n$, then $R$ lies in $\hat{C}$. 
\end{itemize}
\end{proposition}

\textbf{Proof.} 
Suppose that $R(t)$, $t \in [0,T)$, is a solution of the ODE $\frac{d}{dt} R(t) = Q(R(t))$ with $R(0) \in \hat{C}$. Then $\hat{R}(t)$, $t \in [0,T)$, is a solution of the analogous ODE on $\mathbb{R}^n \times \mathbb{R}^2$. Since $\hat{R}(0)$ has nonnegative isotropic curvature, Proposition \ref{pic.is.preserved.2} implies that $\hat{R}(t)$ has nonnegative isotropic curvature for all $t \in [0,T)$. Thus, we conclude that $R(t) \in \hat{C}$ for all $t \in [0,T)$. 

In order to prove (ii), we consider an algebraic curvature operator $R \in \hat{C}$. Let $\{e_1,e_2\}$ be an orthonormal two-frame in $\mathbb{R}^n$. We define an orthonormal four-frame $\{\hat{e}_1,\hat{e}_2,\hat{e}_3,\hat{e}_4\}$ in $\mathbb{R}^n \times \mathbb{R}^2$ by 
\[\begin{array}{l@{\qquad\qquad}l} 
\hat{e}_1 = (e_1,0,0), & \hat{e}_2 = (0,0,1), \\ 
\hat{e}_3 = (e_2,0,0), & \hat{e}_4 = (0,1,0). 
\end{array}\] 
Since $\hat{R}$ has nonnegative isotropic curvature, we have 
\begin{align*} 
0 &\leq \hat{R}(\hat{e}_1,\hat{e}_3,\hat{e}_1,\hat{e}_3) + \hat{R}(\hat{e}_1,\hat{e}_4,\hat{e}_1,\hat{e}_4) \\ 
&+ \hat{R}(\hat{e}_2,\hat{e}_3,\hat{e}_2,\hat{e}_3) + \hat{R}(\hat{e}_2,\hat{e}_4,\hat{e}_2,\hat{e}_4) - 2 \, \hat{R}(\hat{e}_1,\hat{e}_2,\hat{e}_3,\hat{e}_4) \\ 
&= R(e_1,e_2,e_1,e_2). 
\end{align*} 
Hence, $R$ has nonnegative sectional curvature. 

It remains to verify (iii). Let $R$ be a nonnegative curvature operator on $\mathbb{R}^n$. Let $\{\hat{e}_1,\hat{e}_2,\hat{e}_3,\hat{e}_4\}$ be an orthonormal four-frame in $\mathbb{R}^n \times \mathbb{R}^2$. We write $\hat{e}_j = (v_j,x_j)$, where $v_j \in \mathbb{R}^n$ and $x_j \in \mathbb{R}^2$. Letting 
\begin{align*} 
\varphi &= v_1 \wedge v_3 + v_4 \wedge v_2, \\ 
\psi &= v_1 \wedge v_4 + v_2 \wedge v_3, 
\end{align*}
we obtain 
\begin{align*} 
&\hat{R}(\hat{e}_1,\hat{e}_3,\hat{e}_1,\hat{e}_3) + \hat{R}(\hat{e}_1,\hat{e}_4,\hat{e}_1,\hat{e}_4) \\ 
&+ \hat{R}(\hat{e}_2,\hat{e}_3,\hat{e}_2,\hat{e}_3) + \hat{R}(\hat{e}_2,\hat{e}_4,\hat{e}_2,\hat{e}_4) - 2 \, \hat{R}(\hat{e}_1,\hat{e}_2,\hat{e}_3,\hat{e}_4) \\ 
&= R(\varphi,\varphi) + R(\psi,\psi) \geq 0. 
\end{align*}
Thus, we conclude that $R \in \hat{C}$. \\

We next apply a technique discovered by C.~B\"ohm and B.~Wilking \cite{Bohm-Wilking}. For each pair of real numbers $a,b$, B\"ohm and Wilking define a linear transformation $\ell_{a,b}$ on the space of algebraic curvature operators by 
\[\ell_{a,b}(R) = R + b\  \text{\rm Ric}_0 \owedge \text{\rm id} + \frac{a}{n} \, \text{\rm scal}\  \text{\rm id} \owedge \text{\rm id}.\] 
Here, $\text{\rm scal}$ and $\text{\rm Ric}_0$ denote the scalar curvature and trace-free Ricci tensor of $R$, respectively. Moreover, $\owedge$ denotes the Kulkarni-Nomizu product, i.e. 
\[(A \owedge B)_{ijkl} = A_{ik} \, B_{jl} - A_{il} \, B_{jk} - A_{jk} \, B_{il} + A_{jl} \, B_{ik}.\] 
For abbreviation, let $I = \frac{1}{2} \: \text{\rm id} \owedge \text{\rm id}$. Combining Proposition \ref{properties.of.C.hat} with results of B\"ohm and Wilking \cite{Bohm-Wilking} yields: 

\begin{proposition} 
\label{pinching.family.1}
Let $0 < b \leq \frac{1}{2}$. We define a cone $\hat{C}(b)$ by 
\[\hat{C}(b) = \{\ell_{a,b}(R): \text{$R \in \hat{C}$ and $\text{\rm Ric} \geq \frac{p}{n} \, \text{\rm scal}$}\},\] 
where \[2a = \frac{2b + (n-2)b^2}{1 + (n-2)b^2}, \qquad p = 1 - \frac{1}{1+(n-2)b^2}.\] Then the cone $\hat{C}(b)$ is transversally invariant under the ODE $\frac{d}{dt} R = Q(R)$. More precisely, for each $R \in \partial \hat{C}(b) \setminus \{0\}$, $Q(R)$ lies in the interior of the tangent cone to $\hat{C}(b)$ at $R$. 
\end{proposition}

\begin{proposition} 
\label{pinching.family.2}
Let $a > \frac{1}{2}$. We define a cone $\hat{C}(a)$ by 
\[\hat{C}(a) = \{\ell_{a,b}(R): \text{$R \in \hat{C}$ and $\text{\rm Ric} \geq \frac{p}{n} \, \text{\rm scal}$}\},\] 
where \[b = \frac{1}{2}, \qquad p = 1 - \frac{4}{n-2+8a}.\] Then the cone $\hat{C}(a)$ is transversally invariant under the ODE $\frac{d}{dt} R = Q(R)$. More precisely, for each $R \in \partial \hat{C}(a) \setminus \{0\}$, $Q(R)$ lies in the interior of the tangent cone to $\hat{C}(a)$ at $R$.
\end{proposition} 

The proofs of Proposition \ref{pinching.family.1} and Proposition \ref{pinching.family.2} are analogous to Lemma 3.4 and Lemma 3.5 in \cite{Bohm-Wilking}, respectively.

Proposition \ref{pinching.family.1} and Proposition \ref{pinching.family.2} provide a continuous family $\hat{C}(s)$, $s > 0$, of closed, convex, $O(n)$-invariant cones. It is easy to see that these cones form a pinching family in the sense of B\"ohm and Wilking \cite{Bohm-Wilking}:

\begin{proposition} 
\label{pinching.family.3}
The cones $\hat{C}(s)$, $s > 0$, have the following properties: 
\begin{itemize}
\item[(i)] For each $R \in \partial \hat{C}(s) \setminus \{0\}$, $Q(R)$ lies in the interior of the tangent cone to $\hat{C}(s)$ at $R$. 
\item[(ii)] $I$ lies in the interior of $\hat{C}(s)$.
\item[(iii)] Given any $\delta \in (0,1)$, there exists a real number $s > 0$ such that every algebraic curvature operator $R \in \hat{C}(s) \setminus \{0\}$ is $\delta$-pinched.
\end{itemize}
\end{proposition}



The convergence of the normalized Ricci flow follows now from a result of B\"ohm and Wilking (cf. \cite{Bohm-Wilking}, Theorem 5.1) which in turn relies on work of Hamilton (cf. \cite{Hamilton2}, Section 5). The proof of that result requires the construction of a suitable pinching set for the ODE. We have a slightly different construction of such a set, which we provide for the convenience of the reader.

\begin{proposition} 
\label{pinching.family.4}
Fix a compact interval $[\alpha,\beta] \subset (0,\infty)$. Assume that $F_0$ is a closed set which is invariant under the ODE $\frac{d}{dt} R = Q(R)$. Moreover, suppose that 
 \[F_0 \subset \{R: R + hI \in \hat{C}(s)\}\] 
for some $s \in [\alpha,\beta]$ and some $h > 0$. Then there exists a real number $\varepsilon > 0$, depending only on $\alpha$, $\beta$, and $n$, such that the following hold: 
\begin{itemize}
\item[(i)] The set  \[F_1 = F_0 \cap \{R: R + 2hI \in \hat{C}(s+\varepsilon)\}\] is invariant under the ODE $\frac{d}{dt} R = Q(R)$.
\item[(ii)] We have \[F_0 \cap \{R: \text{\rm tr}(R) \leq h\} \subset F_1.\]
\end{itemize}
\end{proposition}

\textbf{Proof.} 
For each $R \in \partial \hat{C}(s) \setminus \{0\}$, $Q(R)$ lies in the interior of the tangent cone to $\hat{C}(s)$ at $R$. Since $Q(R)$ is homogenous of degree $2$, we can find a constant $N > 1$, depending only on $\alpha$, $\beta$, and $n$, with the following property: if $R \in \partial \hat{C}(s)$ for some $s \in [\alpha,\beta+1]$ and $\text{\rm tr}(R) > N$, then $Q(R - 2I)$ lies in the interior of the tangent cone to $\hat{C}(s)$ at $R$. 

Observe that $\{R: R + I \in \hat{C}(s)\} \cap \{R: \text{\rm tr}(R) \leq N\}$ is a compact set which is contained in the interior of the set $\{R: R + 2I \in \hat{C}(s)\}$. Hence, there exists a real number $\varepsilon \in (0,1)$, depending only on $\alpha$, $\beta$, and $n$, such that 
 \[\{R: R + I \in \hat{C}(s)\} \cap \{R: \text{\rm tr}(R) \leq N\} \subset \{R: R + 2I \in \hat{C}(s+\varepsilon)\}\]
for all $s \in [\alpha,\beta]$. 

We now consider the set 
\[F_1 = F_0 \cap \{R: R + 2hI \in \hat{C}(s+\varepsilon)\}.\] 
Using the inclusions 
\[F_0 \subset \{R: R + hI \in \hat{C}(s)\}\] 
and
\[\{R: R + hI \in \hat{C}(s)\} \cap \{R: \text{\rm tr}(R) \leq Nh\} \subset \{R: R + 2hI \in \hat{C}(s+\varepsilon)\},\] 
we obtain 
\[F_0 \cap \{R: \text{\rm tr}(R) \leq Nh\} \subset F_1.\] 
Hence, it remains to show that the set $F_1$ is invariant under the ODE $\frac{d}{dt} R = Q(R)$. Let $R(t)$, $t \in [0,T)$, be a solution of the ODE $\frac{d}{dt} R(t) = Q(R(t))$ with $R(0) \in F_1$. Since $F_0$ is invariant under the ODE $\frac{d}{dt} R = Q(R)$, we have $R(t) \in F_0$ for all $t \in [0,T)$. We claim that $R(t) + 2hI \in \hat{C}(s+\varepsilon)$ for all $t \in [0,T)$. Suppose this is false. We then define 
\[\tau = \inf \{t \in [0,T): R(t) + 2hI \notin \hat{C}(s+\varepsilon)\}.\] 
Clearly, $R(\tau) + 2hI \in \partial \hat{C}(s+\varepsilon)$. Moreover, we have $\text{\rm tr}(R(\tau)) \geq Nh$, hence $\text{\rm tr}(R(\tau) + 2hI) > Nh$. Consequently, $Q(R(\tau))$ lies in the interior of the tangent cone to $\hat{C}(s+\varepsilon)$ at $R(\tau) + 2hI$. This contradicts the definition of $\tau$. Thus, we conclude that $R(t) \in F_1$ for all $t \in [0,T)$. \\

\begin{proposition}[B\"ohm and Wilking \cite{Bohm-Wilking}, Theorem 4.1]
\label{pinching.family.5}
Suppose that $K$ is a compact set which is contained in the interior of $\hat{C}$. Then there exists a closed, convex, $O(n)$-invariant set $F$ with the following properties: 
\begin{itemize}
\item[(i)] $F$ is invariant under the ODE $\frac{d}{dt} R = Q(R)$. 
\item[(ii)] For each $\delta \in (0,1)$, the set $\{R \in F: \text{$R$ is not $\delta$-pinched}\}$ is bounded.
\item[(iii)] $K$ is a subset of $F$. 
\end{itemize}
\end{proposition}

\textbf{Proof.} 
By scaling, we may assume that $\text{\rm tr}(R) \leq 1$ for all $R \in K$. Since $K$ is contained in the interior of $\hat{C}$, there exists a positive real number $s_0$ such that $K \subset \hat{C}(s_0)$. We now apply Proposition \ref{pinching.family.4} with $F_0 = \hat{C}(s_0)$ and $h = 1$. Hence, there exists a real number $s_1 > s_0$ such that the set  \[F_1 = F_0 \cap \{R: R + 2I \in \hat{C}(s_1)\}\] is invariant under the ODE $\frac{d}{dt} R = Q(R)$, and 
\[F_0 \cap \{R: \text{\rm tr}(R) \leq 1\} \subset F_1.\] 
Proceeding inductively, we obtain an increasing sequence of real numbers $s_j$, $j \in \mathbb{N}$, and a sequence of closed, convex, $O(n)$-invariant sets $F_j$, $j \in \mathbb{N}$, with the following properties:  
\begin{itemize}
\item[(a)] For each $j \in \mathbb{N}$, we have $F_{j+1} = F_j \cap \{R: R + 2^{j+1} I \in \hat{C}(s_{j+1})\}$. 
\item[(b)] For each $j \in \mathbb{N}$, we have $F_j \cap \{R: \text{\rm tr}(R) \leq 2^j\} \subset F_{j+1}$. 
\item[(c)] For each $j \in \mathbb{N}$, the set $F_j$ is invariant under the ODE $\frac{d}{dt} R = Q(R)$. 
\item[(d)] $s_j \to \infty$ as $j \to \infty$.
\end{itemize}
We now define $F = \bigcap_{j=1}^\infty F_j$. Clearly, $F$ is a closed, convex, $O(n)$-invariant set, which is invariant under the ODE $\frac{d}{dt} R = Q(R)$. Since $K \subset F_0 \cap \{R: \text{\rm tr}(R) \leq 1\}$, it follows from property (b) that $K \subset F_j$ for all $j \in \mathbb{N}$. Hence, $K$ is a subset of $F$. Finally, property (a) implies \[F \subset F_j \subset \{R: R + 2^j I \in \hat{C}(s_j)\}\] for all $j \in \mathbb{N}$. Since $s_j \to \infty$ as $j \to \infty$, the assertion follows from Proposition \ref{pinching.family.3}. \\

Having established the existence of suitable pinching sets, the convergence of the flow follows from the same arguments as in \cite{Bohm-Wilking}, \cite{Hamilton2}:

\begin{theorem} 
\label{main.convergence.theorem}
Let $(M,g_0)$ be a compact Riemannian manifold of dimension $n \geq 4$. Assume that the curvature tensor of $(M,g_0)$ lies in the interior of the cone $\hat{C}$ for all points in $M$. Then the normalized Ricci flow with initial metric $g_0$ exists for all time and converges to a metric of constant sectional curvature as $t \to \infty$.
\end{theorem}

\section{An algebraic characterization of the cone $\hat{C}$}

In this section, we provide a necessary and sufficient condition for $\hat{R}$ to have nonnegative isotropic curvature. We will need the following linear algebra result (cf. \cite{Chen}, Lemma 3.1). We give a short proof of this for completeness.

\begin{lemma} 
\label{linear.algebra.lemma.1}
Assume that $\varphi,\psi \in \wedge^2 \mathbb{R}^4$ are two-vectors satisfying $\varphi \wedge \varphi = \psi \wedge \psi$, $\varphi \wedge \psi = 0$, and $\langle \varphi,\psi \rangle = 0$. Then there exists an orthonormal basis $\{e_1,e_2,e_3,e_4\}$ of $\mathbb{R}^4$ such that 
\begin{align*} 
\varphi &= a_1 \, e_1 \wedge e_3 + a_2 \, e_4 \wedge e_2, \\ 
\psi &= b_1 \, e_1 \wedge e_4 + b_2 \, e_2 \wedge e_3 
\end{align*} 
with $a_1a_2 = b_1b_2$.
\end{lemma}

\textbf{Proof.}
We first consider the (generic) case in which at least one of $\varphi$, $\psi$ is neither self-dual nor anti-self-dual. Without loss of generality, we may assume that $\varphi$ is neither self-dual
nor anti-self-dual. Consider the anti-symmetric bilinear form defined on $\mathbb{R}^4$ 
by $(v,w) \mapsto \langle \varphi,v \wedge w \rangle$. A standard result in linear algebra implies that there exists a positively oriented orthonormal basis $\{v_1,v_2,v_3,v_4\}$ in which $\varphi$ has the form
\[\varphi = a_1 \, v_1 \wedge v_3 + a_2 \, v_4 \wedge v_2\] 
for suitable coefficients $a_1,a_2$. By assumption, we have $\langle \varphi,\psi \rangle = 0$ and $\varphi \wedge \psi = 0$. This implies 
\[a_1 \, \langle \psi,v_1 \wedge v_3 \rangle + a_2 \, \langle \psi,v_4 \wedge v_2 \rangle = 0\] 
and 
\[a_2 \, \langle \psi,v_1 \wedge v_3 \rangle + a_1 \, \langle \psi,v_4 \wedge v_2 \rangle = 0.\] 
Since $\varphi$ is neither self-dual nor anti-self-dual, we have $a_1^2 \neq a_2^2$. Therefore, we obtain $\langle \psi,v_1 \wedge v_3 \rangle = \langle \psi,v_4 \wedge v_2 \rangle = 0$.

We now consider the two-dimensional subspaces $W,Z \subset \mathbb{R}^4$ where $W$ is the span 
of $\{v_1,v_3\}$ and $Z$ is the span of $\{v_4,v_2\}$. We take the orientations on these spaces 
so that the indicated bases are positively oriented. We consider the bilinear pairing $\sigma: W \times Z \to \mathbb{R}$ given by $\sigma(w,z) = \langle \psi,w \wedge z \rangle$. Linear algebra (singular value decomposition) allows us to find 
positively oriented orthonormal bases $\{e_1,e_3\}$ for $W$ and $\{e_4,e_2\}$ for $Z$ such that $\sigma(e_1,e_2) = 0$ and 
$\sigma(e_3,e_4) = 0$. Clearly, $\{e_1,e_2,e_3,e_4\}$ is a positively oriented orthonormal basis of 
$\mathbb{R}^4$. Since $e_1 \wedge e_3 = v_1 \wedge v_3$ and $e_4 \wedge e_2 = v_4 \wedge v_2$, we have 
\[\varphi = a_1 \, e_1 \wedge e_3 + a_2 \, e_4 \wedge e_2.\] 
Moreover, we have 
\begin{align*} 
&\langle \psi,e_1 \wedge e_3 \rangle = \langle \psi,v_1 \wedge v_3 \rangle = 0, \\ 
&\langle \psi,e_4 \wedge e_2 \rangle = \langle \psi,v_4 \wedge v_2 \rangle = 0 
\end{align*} 
and 
\begin{align*} 
&\langle \psi,e_1 \wedge e_2 \rangle = \sigma(e_1,e_2) = 0, \\ 
&\langle \psi,e_3 \wedge e_4 \rangle = \sigma(e_3,e_4) = 0. 
\end{align*} 
Thus, we conclude that 
\[\psi = b_1 \, e_1 \wedge e_4 + b_2 \, e_2 \wedge e_3\] 
for suitable coefficients $b_1,b_2$. The condition $\varphi \wedge \varphi = \psi \wedge \psi$ then implies $a_1a_2 = b_1b_2$.

We next consider the case in which each of $\varphi$ and $\psi$ is either self-dual or anti-self-dual. The condition 
$\varphi \wedge \varphi = \psi \wedge \psi$ implies that they are either both self-dual or both 
anti-self-dual. Without loss of generality assume both are self-dual. Since the assertion is trivial 
for $\varphi = \psi = 0$, we may assume that $\varphi \neq 0$. As above, we choose a positively oriented orthonormal basis $\{v_1,v_2,v_3,v_4\}$ in which $\varphi = a(v_1 \wedge v_3 + v_4 \wedge v_2)$ for some $a \neq 0$. The condition $\langle \varphi,\psi \rangle = 0$ implies $\langle \psi,v_1 \wedge v_3 + v_4 \wedge v_2 \rangle = 0$. Since $\psi$ is self-dual, it follows that $\langle \psi,v_1 \wedge v_3 \rangle = \langle \psi,v_4 \wedge v_2 \rangle = 0$. Therefore, we can complete the argument as above. This finishes the proof. \\

\begin{lemma}
\label{linear.algebra.lemma.2}
Assume that $\varphi,\psi \in \wedge^2\mathbb{R}^4$ are two-vectors satisfying $\varphi \wedge \varphi = \psi \wedge \psi$ and 
$\varphi \wedge \psi = 0$. Then there exists an orthonormal basis $\{e_1,e_2,e_3,e_4\}$ of $\mathbb{R}^4$ 
and real numbers $a_1,a_2,b_1,b_2,\theta$ such that $a_1a_2 = b_1b_2$ and 
\begin{align*} 
\cos \theta \, \varphi + \sin \theta \, \psi &= a_1 \, e_1 \wedge e_3 + a_2 \, e_4 \wedge e_2, \\ 
-\sin \theta \, \varphi + \cos \theta \, \psi &= b_1 \, e_1 \wedge e_4 + b_2 \, e_2 \wedge e_3. 
\end{align*} 
\end{lemma}

\textbf{Proof.} 
We choose a real number $\theta$ such that \[\frac{1}{2} \sin(2\theta) \, (|\varphi|^2 - |\psi|^2) = \cos(2\theta) \, \langle \varphi,\psi \rangle.\] 
We then define 
\begin{align*} 
\tilde{\varphi} &= \cos \theta \, \varphi + \sin \theta \, \psi, \\ 
\tilde{\psi} &= -\sin \theta \, \varphi + \cos \theta \, \psi. 
\end{align*} 
By assumption, we have $\varphi \wedge \varphi = \psi \wedge \psi$ and $\varphi \wedge \psi = 0$. This implies 
\[\tilde{\varphi} \wedge \tilde{\varphi} - \tilde{\psi} \wedge \tilde{\psi} 
= \cos(2\theta) \, (\varphi \wedge \varphi - \psi \wedge \psi) + 2\sin(2\theta) \, \varphi \wedge \psi = 0\] 
and 
\[\tilde{\varphi} \wedge \tilde{\psi} = -\frac{1}{2} \sin(2\theta) \, 
(\varphi \wedge \varphi - \psi \wedge \psi) + \cos(2\theta) \, \varphi \wedge \psi = 0.\] 
Moreover, we have 
\[\langle \tilde{\varphi},\tilde{\psi} \rangle = -\frac{1}{2} \sin(2\theta) \, 
(|\varphi|^2 - |\psi|^2) + \cos(2\theta) \, \langle \varphi,\psi \rangle = 0\] 
by definition of $\theta$. Hence, the assertion follows from Lemma \ref{linear.algebra.lemma.1}. \\

\begin{proposition} 
\label{characterization.of.C.hat}
Let $R$ be an algebraic curvature operator on $\mathbb{R}^n$, and let $\hat{R}$ be the induced algebraic curvature operator on $\mathbb{R}^n \times \mathbb{R}^2$. The following statements are equivalent: 
\begin{itemize}
\item[(i)] $\hat{R}$ has nonnegative isotropic curvature.
\item[(ii)] For all orthonormal four-frames $\{e_1,e_2,e_3,e_4\}$ and all $\lambda,\mu \in [-1,1]$, we have 
\begin{align*} 
&R(e_1,e_3,e_1,e_3) + \lambda^2 \, R(e_1,e_4,e_1,e_4) \\ 
&+ \mu^2 \, R(e_2,e_3,e_2,e_3) + \lambda^2\mu^2 \, R(e_2,e_4,e_2,e_4) - 2\lambda\mu \, R(e_1,e_2,e_3,e_4) \geq 0. 
\end{align*} 
\end{itemize}
\end{proposition}

\textbf{Proof.} 
Assume first that $\hat{R}$ has nonnegative isotropic curvature. Let 
$\{e_1,e_2,e_3,e_4\}$ be an orthonormal four-frame in $\mathbb{R}^n$, and let $\lambda,\mu \in [-1,1]$. We define 
\[\begin{array}{l@{\qquad\qquad}l} 
\hat{e}_1 = (e_1,0,0), & \hat{e}_2 = (\mu e_2,0,\sqrt{1-\mu^2}), \\ 
\hat{e}_3 = (e_3,0,0), & \hat{e}_4 = (\lambda e_4,\sqrt{1-\lambda^2},0). 
\end{array}\] 
Clearly, the vectors $\{\hat{e}_1,\hat{e}_2,\hat{e}_3,\hat{e}_4\}$ form an orthonormal four-frame in 
$\mathbb{R}^n \times \mathbb{R}^2$. Since $\hat{R}$ has nonnegative isotropic curvature, we have 
\begin{align*} 
0 &\leq \hat{R}(\hat{e}_1,\hat{e}_3,\hat{e}_1,\hat{e}_3) + \hat{R}(\hat{e}_1,\hat{e}_4,\hat{e}_1,\hat{e}_4) \\ 
&+ \hat{R}(\hat{e}_2,\hat{e}_3,\hat{e}_2,\hat{e}_3) + \hat{R}(\hat{e}_2,\hat{e}_4,\hat{e}_2,\hat{e}_4) - 2 \, \hat{R}(\hat{e}_1,\hat{e}_2,\hat{e}_3,\hat{e}_4) \\ 
&= R(e_1,e_3,e_1,e_3) + \lambda^2 \, R(e_1,e_4,e_1,e_4) \\ 
&+ \mu^2 \, R(e_2,e_3,e_2,e_3) + \lambda^2\mu^2 \, R(e_2,e_4,e_2,e_4) - 2\lambda\mu \, R(e_1,e_2,e_3,e_4), 
\end{align*} 
as claimed.

Conversely, assume that (ii) holds. We claim that $\hat{R}$ has nonnegative isotropic curvature. 
Let $\{\hat{e}_1,\hat{e}_2,\hat{e}_3,\hat{e}_4\}$ be an orthonormal four-frame in $\mathbb{R}^n 
\times \mathbb{R}^2$. We write $\hat{e}_j = (v_j,x_j)$, where $v_j \in \mathbb{R}^n$ and $x_j \in \mathbb{R}^2$. 
Let $V$ be a four-dimensional subspace of $\mathbb{R}^n$ containing $\{v_1,v_2,v_3,v_4\}$. We define 
\begin{align*} 
\varphi &= v_1 \wedge v_3 + v_4 \wedge v_2 \in \wedge^2 V, \\ 
\psi &= v_1 \wedge v_4 + v_2 \wedge v_3 \in \wedge^2 V. 
\end{align*}
Clearly, $\varphi \wedge \varphi = \psi \wedge \psi$ and $\varphi \wedge \psi = 0$. By Lemma \ref{linear.algebra.lemma.2}, 
there exists an orthonormal basis $\{e_1,e_2,e_3,e_4\}$ of $V$ and real numbers $a_1,a_2,b_1,b_2,\theta$ such that $a_1a_2 = b_1b_2$ and 
\begin{align*} 
\tilde{\varphi} &:= \cos \theta \, \varphi + \sin \theta \, \psi = a_1 \, e_1 \wedge e_3 + a_2 \, e_4 \wedge e_2, \\ 
\tilde{\psi} &:= -\sin \theta \, \varphi + \cos \theta \, \psi = b_1 \, e_1 \wedge e_4 + b_2 \, e_2 \wedge e_3. 
\end{align*} 
Using the first Bianchi identity, we obtain 
\begin{align*} 
R(\varphi,\varphi) + R(\psi,\psi) 
&= R(\tilde{\varphi},\tilde{\varphi}) + R(\tilde{\psi},\tilde{\psi}) \\ 
&= a_1^2 \, R(e_1,e_3,e_1,e_3) + b_1^2 \, R(e_1,e_4,e_1,e_4) \\ 
&+ b_2^2 \, R(e_2,e_3,e_2,e_3) + a_2^2 \, R(e_2,e_4,e_2,e_4) \\ 
&- 2a_1a_2 \, R(e_1,e_2,e_3,e_4). 
\end{align*} 
Condition (ii) implies that the right hand side is nonnegative. From this, it follows that 
\begin{align*} 
&\hat{R}(\hat{e}_1,\hat{e}_3,\hat{e}_1,\hat{e}_3) + \hat{R}(\hat{e}_1,\hat{e}_4,\hat{e}_1,\hat{e}_4) \\ 
&+ \hat{R}(\hat{e}_2,\hat{e}_3,\hat{e}_2,\hat{e}_3) + \hat{R}(\hat{e}_2,\hat{e}_4,\hat{e}_2,\hat{e}_4) - 2 \, \hat{R}(\hat{e}_1,\hat{e}_2,\hat{e}_3,\hat{e}_4) \\ 
&= R(\varphi,\varphi) + R(\psi,\psi) \geq 0. 
\end{align*}
Hence, $\hat{R}$ has nonnegative isotropic curvature. \\

\begin{corollary} 
\label{pinching}
Assume that all sectional curvatures of $R$ lie in the interval $[1,4]$. Then $\hat{R}$ has nonnegative isotropic curvature. 
\end{corollary}

\textbf{Proof.} 
Let $\{e_1,e_2,e_3,e_4\}$ be an orthonormal four-frame, and let $\lambda,\mu \in [-1,1]$. Since the 
sectional curvatures of $R$ lie in the interval $[1,4]$, we have $|R(e_1,e_2,e_3,e_4)| \leq 2$ 
by Berger's inequality (see e.g. \cite{Karcher}). Thus, we conclude that 
\begin{align*} 
&R(e_1,e_3,e_1,e_3) + \lambda^2 \, R(e_1,e_4,e_1,e_4) \\ 
&+ \mu^2 \, R(e_2,e_3,e_2,e_3) + \lambda^2\mu^2 \, R(e_2,e_4,e_2,e_4) - 2\lambda\mu \, R(e_1,e_2,e_3,e_4) \\ 
&\geq 1 + \lambda^2 + \mu^2 + \lambda^2 \mu^2 - 4 \, |\lambda\mu| \\ 
&\geq 0. 
\end{align*} 
Hence, the assertion follows from Proposition \ref{characterization.of.C.hat}.

\end{document}